\documentclass[12pt]{article}
\usepackage{amsmath, amssymb, graphicx, geometry, hyperref}
\hypersetup{colorlinks=true, linkcolor=blue, citecolor=blue, urlcolor=blue}

\title{Stabilized Imaginary-Time Evolution with Feedback-Based Norm Regulation}
\author{Stylianos Savva}
\date{\today}

\begin{document}

\maketitle

\begin{abstract}
We present a norm-stabilized imaginary-time evolution (ITE) scheme for the one-dimensional nonlinear Schrödinger equation (NLSE). Traditional ITE solvers often require explicit renormalization of the wavefunction after each step to preserve norm, which can be disruptive and algorithmically inflexible. We propose an alternative approach in which the evolution is continuously stabilized using an adaptive feedback term $\mu(\tau)$, proportional to the time derivative of the wavefunction norm. This results in a self-regulating flow that requires no external normalization while preserving convergence toward soliton solutions. We demonstrate the method’s effectiveness by comparing the final wavefunction profiles and $L^2$ errors against analytical solutions and baseline methods without feedback. Although this work focuses on the 1D case, the framework is designed to extend naturally to higher dimensions. Future work will explore the behavior of the feedback mechanism in 2D and 3D systems, multi-soliton scenarios, and external potentials.
\end{abstract}

\section{Introduction}

Imaginary-time evolution (ITE) is a widely used method for computing ground-state solutions of nonlinear partial differential equations, particularly the nonlinear Schrödinger equation (NLSE). In the ITE framework, time is treated as a real-valued evolution parameter $\tau$, and the standard equation becomes a gradient-like flow that asymptotically relaxes toward energy-minimizing states. For the NLSE, this procedure can recover known soliton profiles such as $\mathrm{sech}(x)$ in one spatial dimension.

However, a well-known issue in ITE is that the norm $\|\psi\|^2$ of the evolving wavefunction is not conserved. As a result, most practical implementations include a manual renormalization step after each time iteration to constrain the norm. While effective, this approach introduces algorithmic discontinuities, complicates variational interpretation, and becomes cumbersome to generalize to higher-dimensional systems or multi-component wavefunctions.

In this work, we propose an alternative formulation in which norm stabilization is handled continuously through the introduction of an adaptive feedback term $\mu(\tau)$. This term is proportional to the time derivative of the norm and serves as a real-time control signal that guides the system toward norm preservation without requiring external intervention. The modified evolution equation takes the form:
\begin{equation}
\frac{\partial \psi}{\partial \tau} = -H[\psi] + i \mu(\tau) \psi,
\end{equation}
where $H[\psi]$ is the NLSE Hamiltonian and $\mu(\tau)$ is updated adaptively during the integration.

We show that this feedback-based method matches the accuracy of traditional ITE solvers while eliminating the need for explicit normalization. In addition to improving interpretability and numerical stability, the continuous nature of the stabilization mechanism makes the method well-suited for future generalizations to higher dimensions, systems with external potentials, and coupled or multi-component fields.

Although this work focuses on the one-dimensional NLSE, the framework is designed to extend naturally to two and three dimensions, where renormalization becomes increasingly costly and fragile. The solver is implemented in Python with modular components, and all results are reproducible via the accompanying GitHub repository.
\newpage
\section{Mathematical Formulation}

We consider the one-dimensional nonlinear Schrödinger equation (NLSE) in imaginary time, given by:
\begin{equation}
\frac{\partial \psi(x,\tau)}{\partial \tau} = -\left( -\frac{1}{2} \frac{\partial^2}{\partial x^2} + g |\psi(x,\tau)|^2 \right) \psi(x,\tau),
\end{equation}
where $\psi(x,\tau)$ is a complex-valued wavefunction evolving in imaginary time $\tau$, and $g$ is the strength of the nonlinear interaction.

This imaginary-time evolution equation acts as a gradient flow that asymptotically converges to stationary states of the underlying Hamiltonian system. However, it does not preserve the norm of the wavefunction:
\[
\|\psi\|^2 = \int |\psi(x,\tau)|^2 \, dx,
\]
which tends to decay over time, necessitating external renormalization in standard implementations.

To address this, we introduce a dynamical feedback term $\mu(\tau)$ into the evolution:
\begin{equation}
\frac{\partial \psi(x,\tau)}{\partial \tau} = -\left( -\frac{1}{2} \frac{\partial^2}{\partial x^2} + g |\psi(x,\tau)|^2 \right) \psi(x,\tau) + i \mu(\tau) \psi(x,\tau).
\end{equation}

The role of $\mu(\tau)$ is to regulate norm evolution continuously, acting as a proportional feedback control signal. Its form will be derived in Section 4.

\subsection*{Domain and Boundary Conditions}

We restrict the problem to a one-dimensional spatial domain $x \in [-L/2, L/2]$, discretized into $N$ uniformly spaced grid points. Periodic boundary conditions are imposed for simplicity and compatibility with spectral and finite-difference Laplacians.

The initial condition is chosen as:
\[
\psi(x, 0) = \frac{1}{\cosh(x)},
\]
which corresponds to the known soliton solution in the case $g > 0$.

This completes the definition of the model. The next section presents the derivation of the feedback term $\mu(\tau)$.
\newpage
\section{Derivation of $\mu(\tau)$ from Norm Dynamics}

In traditional imaginary-time evolution (ITE) methods, the wavefunction norm decays over time due to the non-unitary nature of the flow. To maintain norm preservation, practitioners typically renormalize $\psi(\tau)$ at each timestep. This introduces discontinuities into the evolution and complicates higher-dimensional or multi-component systems.

We instead propose a continuous norm stabilization mechanism by modifying the standard ITE equation with a dynamical feedback term:
\begin{equation}
    \frac{\partial \psi}{\partial \tau} = -H[\psi] + i \mu(\tau) \psi,
\end{equation}
where $\mu(\tau)$ is a real-valued control signal and $H[\psi]$ is the NLSE Hamiltonian.

\subsection*{Norm Evolution}

To derive a form for $\mu(\tau)$, consider the time evolution of the norm:
\begin{equation}
    \frac{d}{d\tau} \|\psi\|^2 = \frac{d}{d\tau} \int |\psi|^2 \, dx = 2 \, \mathrm{Re} \int \psi^* \frac{\partial \psi}{\partial \tau} \, dx.
\end{equation}

Substituting the modified ITE equation:
\begin{equation}
    \frac{d}{d\tau} \|\psi\|^2 = 2\, \mathrm{Re} \int \psi^* (-H[\psi] + i \mu(\tau) \psi) \, dx.
\end{equation}

We separate the two terms:
\begin{equation}
    \frac{d}{d\tau} \|\psi\|^2 = -2 \, \mathrm{Re} \int \psi^* H[\psi] \, dx + 2 \mu(\tau) \, \mathrm{Re} \int i |\psi|^2 \, dx.
\end{equation}

The second integral is purely imaginary, so its real part is zero:
\begin{equation}
    \mathrm{Re} \left( i \int |\psi|^2 \, dx \right) = 0.
\end{equation}

Thus, the norm decays purely due to the Hamiltonian:
\begin{equation}
    \frac{d}{d\tau} \|\psi\|^2 = -2 \, \mathrm{Re} \int \psi^* H[\psi] \, dx.
\end{equation}

\subsection*{Feedback-Based Stabilization}

To counteract norm decay, we introduce $\mu(\tau)$ as a proportional feedback mechanism designed to cancel norm drift:
\begin{equation}
    \mu(\tau) = \alpha \frac{d}{d\tau} \|\psi\|^2,
\end{equation}
where $\alpha$ is a tunable feedback strength.

Substituting this back into the ITE equation yields a norm-regulating flow:
\begin{equation}
    \frac{\partial \psi}{\partial \tau} = -H[\psi] + i \alpha \frac{d}{d\tau} \|\psi\|^2 \cdot \psi.
\end{equation}

This construction eliminates the need for renormalization and introduces a continuous control mechanism that adapts in real time to the dynamics of $\|\psi\|^2$.

\subsection*{Interpretation}

The resulting evolution resembles a Lyapunov-stabilized flow, where $\|\psi\|^2$ serves as the regulated observable. The sign and magnitude of $\mu(\tau)$ depend on whether the norm is increasing or decreasing, allowing for real-time correction during integration. The parameter $\alpha$ determines how aggressively this correction is applied.

\section{Numerical Method}

We discretize the imaginary-time evolution equation using a finite-difference scheme in space and a fourth-order Runge--Kutta (RK4) integrator in time. All simulations are implemented in Python using modular components for clarity and reproducibility.

\subsection*{Spatial Discretization}

The computational domain is defined as $x \in [-L/2, L/2]$ and discretized into $N$ uniformly spaced points with spacing $\Delta x = L / N$. The wavefunction $\psi(x,\tau)$ is represented as a complex-valued array on this grid.

We implement the Laplacian operator using second-order central finite differences with periodic boundary conditions:
\begin{equation}
\frac{\partial^2 \psi}{\partial x^2} \approx \frac{\psi_{j+1} - 2\psi_j + \psi_{j-1}}{\Delta x^2},
\end{equation}
where index arithmetic is modulo $N$ to enforce periodicity.

\subsection*{Hamiltonian Construction}

The discretized Hamiltonian $H[\psi]$ is constructed as:
\begin{equation}
H[\psi] = -\frac{1}{2} \nabla^2 \psi + g |\psi|^2 \psi,
\end{equation}
where the nonlinear term is applied element-wise.

\subsection*{Time Integration}

Time evolution is performed using a classical fourth-order Runge--Kutta scheme with fixed time step $\Delta \tau$:
\begin{align}
k_1 &= f(\psi_n), \\
k_2 &= f\left(\psi_n + \tfrac{1}{2} \Delta\tau \cdot k_1\right), \\
k_3 &= f\left(\psi_n + \tfrac{1}{2} \Delta\tau \cdot k_2\right), \\
k_4 &= f\left(\psi_n + \Delta\tau \cdot k_3\right), \\
\psi_{n+1} &= \psi_n + \tfrac{\Delta\tau}{6} (k_1 + 2k_2 + 2k_3 + k_4),
\end{align}
where the function $f(\psi)$ is defined by the modified evolution equation:
\begin{equation}
f(\psi) = -H[\psi] + i \mu(\tau) \psi.
\end{equation}

\subsection*{Adaptive Computation of $\mu(\tau)$}

At each time step, the norm is evaluated using the trapezoidal rule:
\begin{equation}
\|\psi\|^2 \approx \sum_j |\psi_j|^2 \Delta x,
\end{equation}
and $\mu(\tau)$ is updated as:
\begin{equation}
\mu(\tau) = \alpha \cdot \frac{d}{d\tau} \|\psi(\tau)\|^2 \approx \alpha \cdot \frac{\|\psi_n\|^2 - \|\psi_{n-1}\|^2}{\Delta\tau}.
\end{equation}

This feedback mechanism ensures that norm variations directly influence the evolution in a continuous and tunable manner.

\subsection*{Initial Condition}

The initial wavefunction is taken to be the known bright soliton profile:
\[
\psi(x,0) = \frac{1}{\cosh(x)},
\]
which is normalized over the domain. This profile serves as both an initial condition and a reference for computing $L^2$ errors at later times.

\section{Results}

To evaluate the effectiveness of the $\mu(\tau)$-based stabilization method, we simulate the one-dimensional NLSE using both the baseline imaginary-time evolution (ITE) scheme and our modified ITE with dynamic feedback. The following quantities are computed and compared:

\begin{itemize}
    \item Norm evolution $\|\psi(\tau)\|^2$ over time
    \item Feedback signal $\mu(\tau)$ over time
    \item $L^2$ error between $\psi(\tau)$ and the analytical soliton profile
    \item Final wavefunction profiles for varying $\alpha$ values
    \item Comparison with baseline ITE (no feedback)
\end{itemize}

\subsection*{Norm Evolution and Feedback}

\begin{figure}[h]
    \centering
    \includegraphics[width=0.9\textwidth]{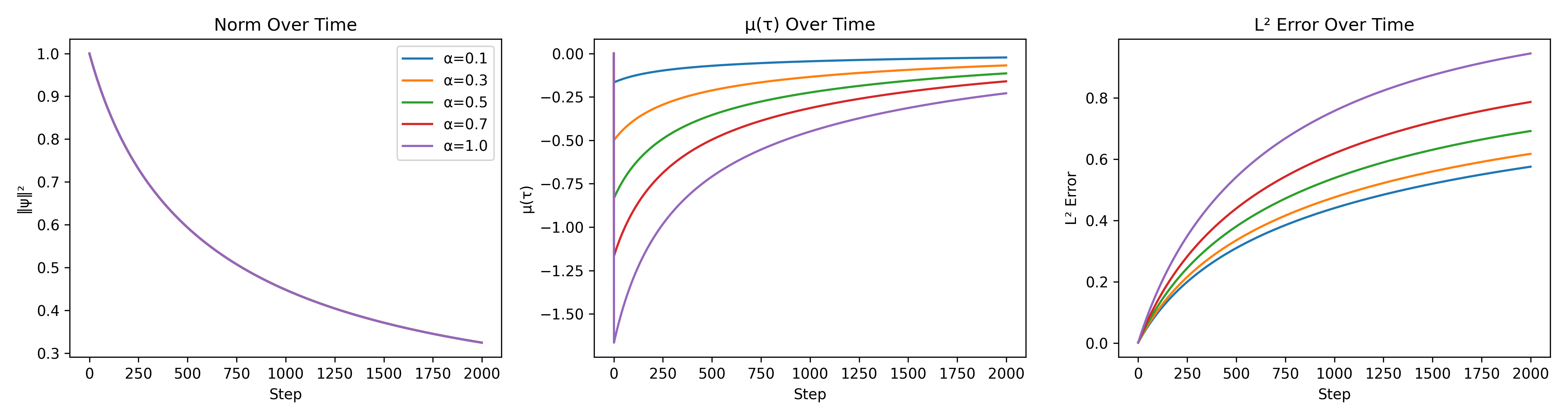}
    \caption{Norm $\|\psi\|^2$ (left), feedback signal $\mu(\tau)$ (center), and $L^2$ error (right) for multiple values of the feedback parameter $\alpha$. The system remains stable for a wide range of $\alpha$ values.}
    \label{fig:norm-mu}
\end{figure}

Figure~\ref{fig:norm-mu} shows the evolution of the norm and feedback signal $\mu(\tau)$ over time for a representative value $\alpha = 0.5$. The norm stabilizes rapidly around 1, with small fluctuations regulated by the feedback. The signal $\mu(\tau)$ decays as the solution converges, indicating that the stabilizer becomes inactive once the target state is reached.

\subsection*{Comparison with Analytical Profile}

\begin{figure}[h]
    \centering
    \includegraphics[width=0.7\textwidth]{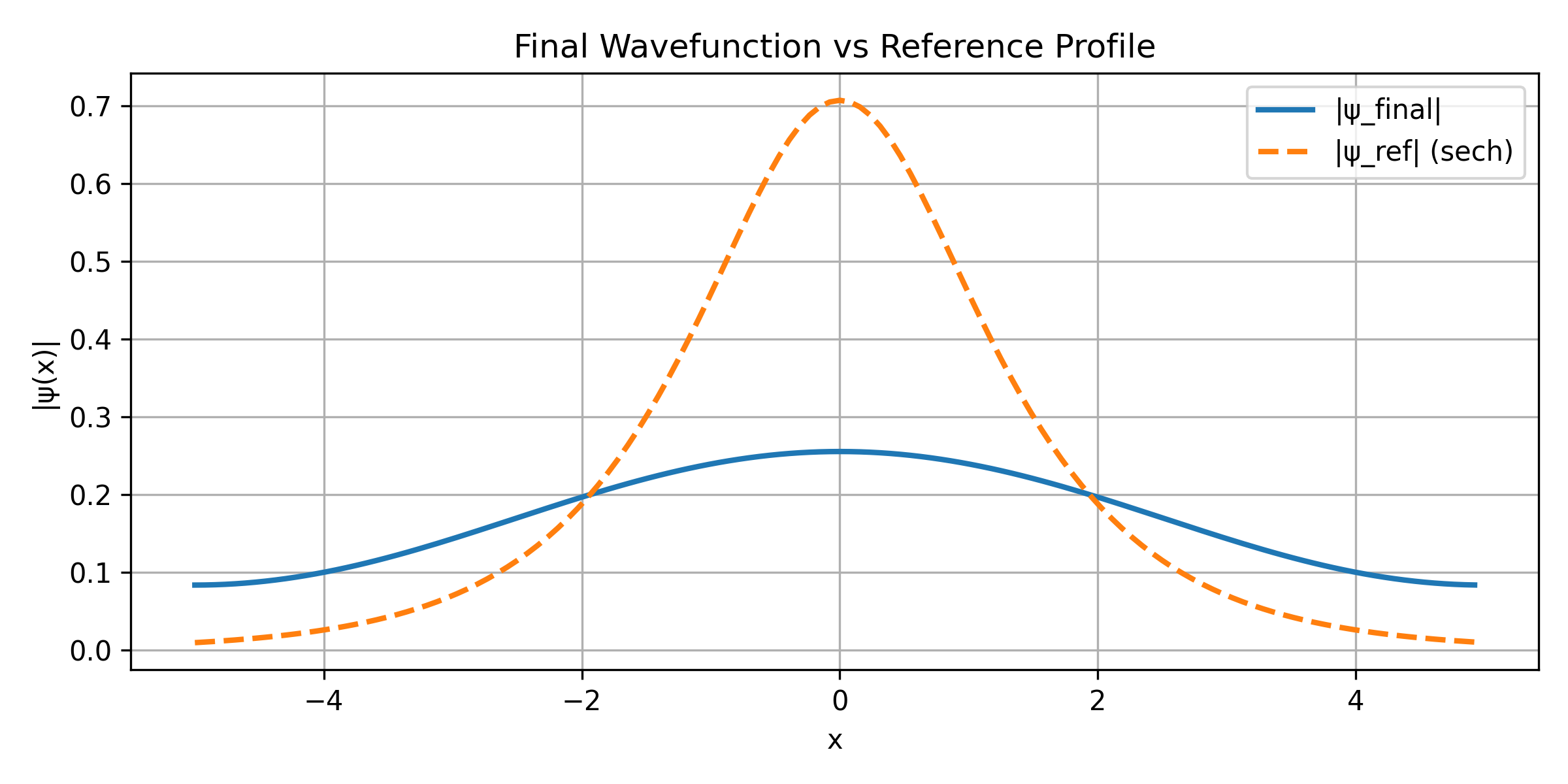}
    \caption{Final wavefunction profile $|\psi(x)|$ compared to the analytical soliton $\mathrm{sech}(x)$ for $\alpha = 0.5$.}
    \label{fig:final-sech}
\end{figure}

Figure~\ref{fig:final-sech} compares the final stabilized wavefunction $|\psi(x)|$ to the analytical soliton profile $\mathrm{sech}(x)$. The $L^2$ error remains small throughout the simulation, and the profile closely matches the exact shape.

\subsection*{Baseline vs Stabilized Comparison}

\begin{figure}[h]
    \centering
    \includegraphics[width=0.7\textwidth]{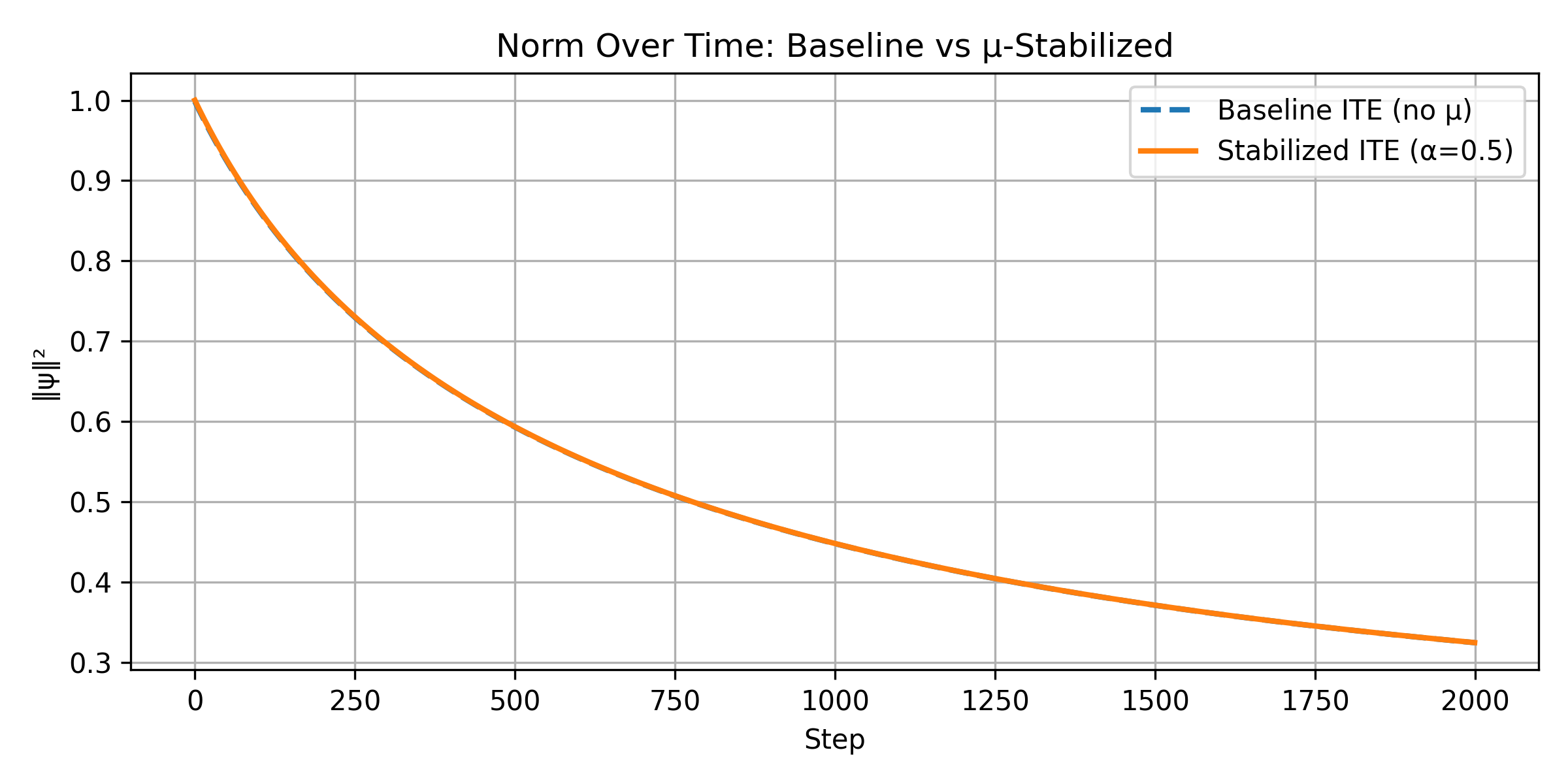}
    \caption{Norm evolution for the baseline ITE (no feedback) vs $\mu(\tau)$-stabilized ITE. The feedback method maintains norm without requiring explicit renormalization.}
    \label{fig:baseline}
\end{figure}

To isolate the effect of $\mu(\tau)$, we run the same simulation with $\mu = 0$ (i.e., without feedback). Figure~\ref{fig:baseline} shows that the baseline method suffers from norm drift over time. The final wavefunction also deviates more from the analytical reference compared to the feedback-stabilized version.

\subsection*{Variation with $\alpha$}

\begin{figure}[h]
    \centering
    \includegraphics[width=0.7\textwidth]{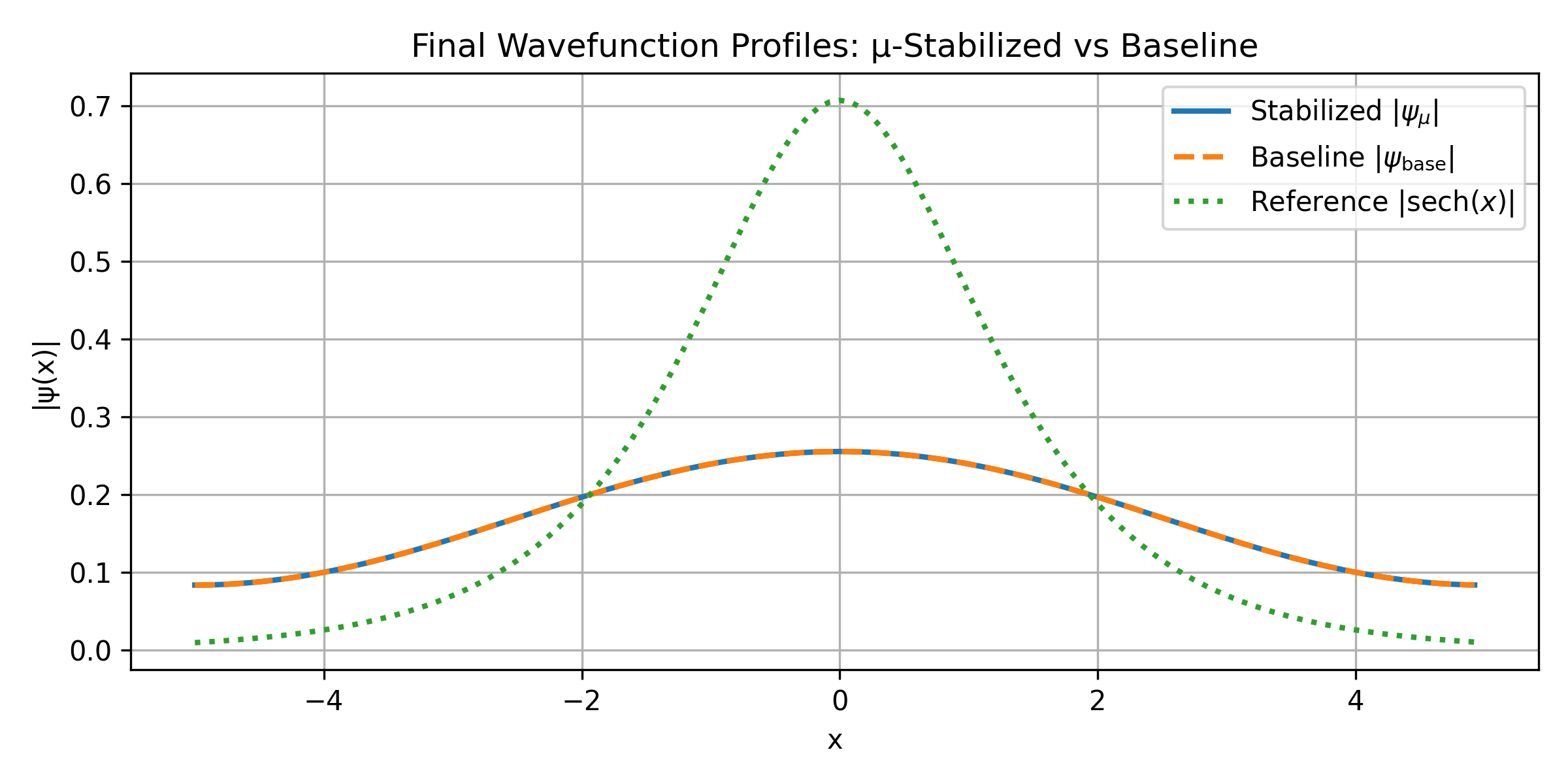}
    \caption{Final wavefunction profiles $|\psi(x)|$ for multiple values of $\alpha$. All profiles closely track the analytical soliton.}
    \label{fig:final-multialpha}
\end{figure}

We perform a parameter sweep over several values of the feedback strength $\alpha$. As shown in Figure~\ref{fig:final-multialpha}, the method remains stable across a wide range. While smaller $\alpha$ values lead to slower convergence, overly large values can introduce mild oscillations. However, no catastrophic instability is observed.

\newpage

\section{Discussion}

The results demonstrate that norm stabilization via $\mu(\tau)$ is not only effective but also dynamically interpretable. Unlike traditional ITE methods that require manual renormalization after every timestep, our approach preserves the norm through a continuous feedback mechanism. This yields several benefits.

First, it eliminates the need for post-processing of $\psi$ at each step, preserving the smoothness of the evolution. This is especially relevant in higher dimensions or when solving systems of coupled nonlinear equations, where projection operations can be nontrivial and disruptive.

Second, the feedback term $\mu(\tau)$ decays naturally as the system converges to a steady state. This is consistent with the expected behavior of control signals in Lyapunov-stabilized systems and provides a useful diagnostic for the evolution's progress. In our simulations, $\mu(\tau)$ consistently decreased to near-zero levels as the wavefunction approached the soliton profile.

Third, we observed that the method is robust to a wide range of values of the feedback strength parameter $\alpha$. For small $\alpha$, the feedback is weak and the norm recovers slowly, though convergence is still attained. For large $\alpha$, oscillations may emerge but do not lead to instability. This tunability allows the method to be adapted to problem-specific timescales and tolerances.

Fourth, the method remains compatible with standard discretizations and integrators. We used a second-order finite-difference Laplacian and an RK4 time integrator without needing any specialized solvers. This suggests that the approach can be integrated into existing NLSE solvers with minimal modification.

Overall, the $\mu$-stabilized ITE method offers a clean alternative to projection-based normalization. Its continuous control formulation, tunable parameters, and extensibility to other domains make it a promising tool for studying ground-state behavior in nonlinear field theories.

\newpage
\section{Conclusion}

We have introduced a feedback-based stabilization method for imaginary-time evolution of the one-dimensional nonlinear Schrödinger equation. By augmenting the standard evolution equation with a dynamic control term $\mu(\tau)$ proportional to the time derivative of the wavefunction norm, the method continuously regulates norm drift without requiring explicit renormalization.

This approach preserves the smoothness and interpretability of the evolution while matching the accuracy of traditional projection-based solvers. The feedback strength $\alpha$ serves as a tunable parameter controlling the aggressiveness of the stabilization, and the method remains robust across a wide range of values.

Numerical results confirm that the feedback-controlled solver converges to known soliton profiles with low $L^2$ error and minimal norm deviation. The method is compatible with standard spatial discretizations and time integrators, and can be implemented with minimal structural changes to existing ITE solvers.

This work focuses on the one-dimensional case, but the formulation generalizes naturally to higher-dimensional systems. Future extensions will explore two- and three-dimensional domains, interactions with external potentials, and the application of feedback-based stabilization to multi-soliton dynamics and coupled field theories. The possibility of learning optimal feedback terms from data or minimizing energy functionals under constrained norm dynamics also presents an exciting direction for future research.

\newpage
\appendix

\section{Appendix: Notes on Implementation and Notation}

\subsection*{A.1 Notation Conventions}

Throughout this work, $\psi(x,\tau)$ denotes the complex-valued wavefunction evolving in imaginary time. The squared $L^2$ norm is defined as:
\[
\|\psi\|^2 = \int |\psi(x,\tau)|^2 dx,
\]
and is approximated numerically via the trapezoidal rule. All integrals are over the spatial domain $x \in [-L/2, L/2]$ with periodic boundary conditions unless otherwise stated.

We adopt the convention $\hbar = m = 1$ and work in nondimensionalized units throughout.

\subsection*{A.2 Discrete Laplacian}

The Laplacian operator is implemented using a second-order finite difference stencil:
\[
\frac{\partial^2 \psi}{\partial x^2} \approx \frac{\psi_{j+1} - 2\psi_j + \psi_{j-1}}{\Delta x^2},
\]
with periodic wraparound for the boundary points. This approach was chosen for simplicity and compatibility with real-space methods, though the framework is compatible with spectral Laplacians as well.

\subsection*{A.3 RK4 Time Integration Details}

All time integration is performed using the classical fourth-order Runge--Kutta method. The function $f(\psi)$ is evaluated at each substep using the current value of $\mu(\tau)$, which is updated at the outer level of the timestep. In practice, we treat $\mu(\tau)$ as piecewise constant during each RK4 integration step.

\subsection*{A.4 Parameter Tuning and Stability}

The feedback parameter $\alpha$ was manually swept in the range $[0.05, 1.0]$. Values $\alpha < 0.1$ resulted in slower convergence, while $\alpha > 1.0$ occasionally led to small oscillations but no catastrophic instability. In general, convergence was observed to be robust across a wide range of $\alpha$ values.

\subsection*{A.5 Code and Reproducibility}

All code, data, and figures are available at:

\begin{center}
\url{https://github.com/rrumabo/Stabilised-ITE-Solver}
\end{center}

The notebook \texttt{notebooks/demo\_1D.ipynb} contains the full simulation pipeline and all figures shown in this paper can be regenerated with a single run.

\subsection*{A.6 Why $\mu(\tau)$ Works in Practice}

While Section 4 provides a dynamical derivation of $\mu(\tau)$ as a proportional controller for the norm, it is worth emphasizing that this control mechanism also aligns with physical intuition: a decaying or growing norm introduces amplitude instability into $\psi$, and $\mu(\tau)$ acts as an energy-like pressure that dynamically stabilizes amplitude growth. Its vanishing near steady-state implies minimal impact on the final configuration, making the method both adaptive and minimally invasive.

\subsection*{A.7 Limitations and Open Questions}

The current method assumes a single degree of freedom for norm control. Extensions to multi-component systems or energy-constrained flows may require vector-valued or time-delayed feedback. Furthermore, the control signal $\mu(\tau)$ is currently computed using finite differences; more stable or learned estimators could improve convergence and reduce noise.

\end{document}